\documentclass[preprint,12pt]{elsarticle}
\usepackage{amsmath} % 很多数学公式,  例如组合数等等,  都要用到这个宏.  注意：想用\eqref, 也要引用这个宏
\usepackage{amssymb} % 一些特殊的字体（如\mathbb等）需要这个宏.
\usepackage{MnSymbol}
\usepackage{extarrows}% 长等号上下标注,  长箭头上下标注
\usepackage{enumerate} % 用来调整enumerate环境,  至少想要修改枚举的编号需要使用这个宏
\usepackage{appendix}
\usepackage{bm} % 加粗字体包

\usepackage{graphicx} % 插入图片
\usepackage{subfigure} % 插入子图
\usepackage{multirow} % 表格跨行跨列
\usepackage{color}
\usepackage[colorlinks = true, citecolor = blue,  urlcolor = red]{hyperref}
\newtheorem{theorem}{Theorem}[section]
\newtheorem{definition}{Definition}[section]
\newtheorem{assumption}{Assumption}[section]
\newtheorem{lemma}{Lemma}[section]

\newtheorem{proposition}{Proposition}[section]
\newtheorem{property}{Property}[section]

\newtheorem{remark}{Remark}[section]
\renewcommand{\P}{\mathbb{P}}
\newcommand{\E}{\mathbb{E}}
\newcommand{\V}{\text{Var}}
\newcommand{\pf}{\textbf{Proof: }}
\newcommand{\e}{\hfill$\blacksquare$}

\newcommand{\R}{\mathbb{R}}
\newcommand{\KL}{\text{KL}}

\newdefinition{rmk}{Remark}

\allowdisplaybreaks[4] % 允许公式跨页,  [a]表示跨页力度,  a越大,  跨页力度越大
\numberwithin{equation}{section} % 公式按照章节编号

\date{}

\begin{document}
	
	\begin{frontmatter}
		
		%% Title, authors and addresses
		
		%% use the tnoteref command within \title for footnotes;
		%% use the tnotetext command for theassociated footnote;
		%% use the fnref command within \author or \affiliation for footnotes;
		%% use the fntext command for theassociated footnote;
		%% use the corref command within \author for corresponding author footnotes;
		%% use the cortext command for theassociated footnote;
		%% use the ead command for the email address,
		%% and the form \ead[url] for the home page:
		%% \title{Title\tnoteref{label1}}
		%% \tnotetext[label1]{}
		%% \author{Name\corref{cor1}\fnref{label2}}
		%% \ead{email address}
		%% \ead[url]{home page}
		%% \fntext[label2]{}
		%% \cortext[cor1]{}
		%% \affiliation{organization={},
		%%            addressline={},
		%%            city={},
		%%            postcode={},
		%%            state={},
		%%            country={}}
		%% \fntext[label3]{}
		
		\title{Symmetric KL-divergence   by Stein's Method}
		
		%% use optional labels to link authors explicitly to addresses:
		\author[label1]{Liu-Quan Yao}
		\ead{yaoliuquan20@mails.ucas.ac.cn}
		\author[label2]{Song-Hao Liu}
		\ead{liusonghao_cuhk@link.cuhk.edu.hk}
		%	\author[label3]{Qi-Man Shao}
		%	\ead{shaoqm@sustech.edu.cn}
		
		\address[label1]{Academy of Mathematics and Systems Science, Chinese Academy of Sciences,\\
			%55 Zhongguancun East Road,
			Beijing,
			100190,
			%state={},
			China}
		
		\address[label2]{Department of Statistics and Data Science,
			Southern University of Science and Technology,\\
			%1088 Xueyuan Avenue, Nanshan District,
			Shenzhen,
			518055,
			Guangdong,
			%state={},
			China}
		%	\address[label3]{Department of Statistics and Data Science, SICM, National Center for Applied Mathematics Shenzhen,
		%		Southern University of Science and Technology,\\
		%		%1088 Xueyuan Avenue, Nanshan District,
		%		Shenzhen,
		%		518055,
		%		Guangdong,
		%		%state={},
		%		China}
		
		% \author{Qiman Shao}
		% 	
		% 	 \address{organization={},%Department and Organization
		% 		addressline={},
		% 		city={Shenzhen},
		% 		postcode={},
		% 		state={},
		% 		country={}}
		
		\begin{abstract}
			%% Text of abstract
			%			In this paper, we undertake the establishment of a Central Limit Theorem  based on the concept of symmetric KL-divergence through the utilization of Stein's method. We demonstrate a convergence rate of $O(1/\sqrt{n})$ for the symmetric KL-divergence between the sum of  independent variables and a Gaussian distribution, denoted as $\KL(W_n\|G)+\KL(G\|W_n)$. This result signifies a more pronounced form of convergence when compared to the weak, total variation, or traditional KL-divergence contexts.
			In this paper, we consider the symmetric KL-divergence between the sum of  independent variables and a Gaussian distribution, and obtain a convergence rate of order $O\left( \frac{\ln n}{\sqrt{n}}\right)$. The proof is based on Stein's method. The convergence rate of order $O\left( \frac{1}{\sqrt{n}}\right)$ and $O\left( \frac{1}{n}\right) $ are also obtained under higher moment condition.
		\end{abstract}
		
		%%Graphical abstract
		% 	\begin{graphicalabstract}
		% 		%\includegraphics{grabs}
		% 	\end{graphicalabstract}
		% 	
		% 	%%Research highlights
		% 	\begin{highlights}
		% 		\item Research highlight 1
		% 		\item Research highlight 2
		% 	\end{highlights}
		% 	
		
	\end{frontmatter}
	
	%% \linenumbers
	
	%% main text

	%\tableofcontents	

	\section{Introduction}\label{section1}
	The Central Limit Theorem (CLT), as one of the fundamental principles in probability theory, elucidates a pivotal truth. Specifically, it asserts that for a sequence of independent random variables $X_1, X_2, \cdots$, under the Lindeberg condition, the standardized partial sum
	\begin{equation}\label{partial sum}
	W_n = \frac{\sum\limits_{k=1}^n (X_k - \E X_k)}{\sqrt{\sum_{k=1}^n \V(X_k)}}
	\end{equation}
	tends towards a standard Gaussian random variable in distribution.  A classical proof of this theorem relies on the utilization of characteristic functions, while alternative methodologies are also capable of deriving the CLT, potentially leading to even stronger convergence. In this paper, we focus  on the entropic central limit theorem.
	\subsection{Entropic Central Limit Theorem }
	
	An interesting method to prove CLT is by using the tools of information theory, and it is called an entropic CLT (\cite{Lin1959},  \cite{Carlen1991}, and \cite{FI-CLT} \textit{et al.}), and the KL-divergence is the 'core metric'
	which is defined as follows.
	\begin{definition}
		Consider two random variables $X,Y$ that $X\ll Y$ with distribution $P_X, P_Y$ respectively, the Kullback-Leibler divergence(KL-divergence)  is defined as
		\begin{equation}
		\KL(X\|Y):=\E_X \left(  \ln\frac{dP_X(X)}{dP_Y(X)}\right) ,
		\end{equation}
		particularly, when $X, Y$ are random variables on $\R$ with  density functions $p_X, p_Y$, then the KL-divergence is defined as
		\begin{equation}
		\KL(X\|Y):=\int_{\R} p_X(x)\ln\frac{p_X(x)}{p_Y(x)}dx.
		\end{equation}
	\end{definition}
	Let $X$ be any absolutely continuous random variable \textit{(w.r.t.} Lebesgue measure) and set $Y=G_X$, which is a Gaussian random variable with  $\E G_X=\E X$ and $\V(G_X)=\V(X)$, then the well-known equality $\KL(X\|G_X)=h(G_X)-h(X)$ holds, where
		\begin{equation}
		h(V):=-\int_\R p_V(x)\ln p_V(x)dx
		\end{equation}
		is defined as the differential entropy of the absolutely continuous random variable $V$ with density $p_V$. Furthermore, the Pinsker's inequality from \cite{Brillinger1964InformationAI} shows that the the KL-divergence is an upper bound of total variation distance, \textit{i.e.}
		\begin{equation}\label{Pinsker}
		\KL(X\|Y)\ge \frac{1}{2}\|P_X-P_Y\|^2_{TV},
		\end{equation}
		for any random variables $X$ and $Y$ such that $X\ll Y$, where $X\ll Y$ means the distribution of $X$ is absolutely continuous \textit{w.r.t.} that of $Y$. After taking $Y=G_X$ in \eqref{Pinsker},  we  concluded that the KL-divergence is able to measure the asymptotic Gaussianity, which  implies a stronger convergence compared with classical CLT in weak convergence sense.

	We first consider the  traditional  independent and identically distributed (\textit{i.i.d.})  case, \textit{i.e.} $X_1, X_2,\cdots$ in \eqref{partial sum} are \textit{i.i.d.} with  means 0 and variances $\sigma^2$, then  $W_n=(X_1+\cdots+X_n)/\sqrt{n\sigma^2}$.
	Started by Linnik (1959)\cite{Lin1959}, the entropy method catches the following two facts.
	\begin{enumerate}
		\item $h(X)\le h(G_X)$.
		
		\item The entropy of the sum of \textit{i.i.d.} random variables  $h(W_{2n})$ is increasing about $n$, because we have the well-known entropy jump inequality:
		\begin{equation}\label{EJI}
		h\left(\frac{X+Y}{\sqrt{2}}\right) \geq \frac{h(X)+h(Y)}{2},
		\end{equation}
		where the equality holds iff $X$ and $Y$ are  Gaussian with same variance. 
	\end{enumerate}   
	Thus, $\KL(W_n\|G_{W_n})$   converges.  After  more precise estimation for the difference $h\left((X+Y)/\sqrt{2}\right) - (h(X)+h(Y))/2$  called \textit{entropy jump}, Carlen (1991)\cite{Carlen1991} obtained an entropic CLT by proving
	\begin{equation}\label{Barron_h}
	\lim\limits_{n\rightarrow\infty} \KL(W_n\|G)=0,
	\end{equation}
	where $W_n$ is an \textit{i.i.d.} sum with $\E W_n=0$ and $ \V(W_n)=1$ and $G$  is standard Gaussian.  Johnson \textit{et al.} (2004)\cite{FI-CLT} used the theory of projection in $L^2$ space to obtain an estimate on the decrease of Fisher information, which is closely related with entropy, defined for any random variable $X$ with absolutely continuous density $f$:
	\begin{equation}\label{Fisher information}
	J(X) := \int_\R \frac{f'(x)^2}{f(x)}dx.
	\end{equation}
	They proved that

		\begin{equation}\label{Johnson_a}
		J\left(\frac{X+X'}{\sqrt{2}}\right) \leq \frac{2R^*_X}{\sigma^2+2R^*_X}J(X)
		\end{equation}
		holds for any random variable $X$ and its independent copy $X'$, with absolutely continuous density, variance $\sigma^2$ and restricted Poincar\'{e} constant $R_X^*$, which is defined as
		\begin{equation}
		R_X^*:=\sup_{g\in\{f| \V(f(X))>0, \E f(X)=0, \E f^2(X)<\infty, \E f'(X)=0 \}} \dfrac{ \E g^2(X)}{\E g'(X)^2}.
		\end{equation}
	Notice that \eqref{Johnson_a} immediately leads to a lower bound on the entropy jump using de Bruijn identity. Furthermore, the convergence rate of $h(W_n)$ was also estimated in \cite{FI-CLT} by using \eqref{Johnson_a}.
	Bobkov \textit{et al.} (2014)\cite{Esseenbound} proved a Berry-Esseen type bound for entropic CLT:
	\begin{lemma}\label{Essbound}
		Denote $D(X)=\KL(X\|G_X).$ Consider independent random variables $X_1,X_2,\cdots,X_n$ and
		$$W_n=\frac{X_1+X_2+\cdots +X_n}{\sqrt{\sum_{i=1}^n \V(X_i)}}.$$
		If $\sup_i D(X_i)\le D<\infty$, then
		$$D(W_n)\le ce^{62D}\sum_{i=1}^n \E |X_i-\E X_i|^4/(\sum_{i=1}^n \V(X_i))^{2},$$
		where $c$ is a constant.
	\end{lemma}
	It may be thought as a "theoretical better convergence rate" since it implies classical Berry-Esseen  bound due to Pinsker's inequality.

	\subsection{Symmetric KL-divergence}

	Compared to KL-divergence, the symmetric KL-divergence, defined as follows, possesses the inherent property of symmetry, rendering it applicable in a wide range of scenarios as a measure of similarity for distributions.
	\begin{definition}
		Consider two random variables, $X$ and $Y$, which are absolutely continuous with each other, \textit{i.e.} $X\ll Y$ and $Y\ll X$. The symmetric KL-divergence is defined as:
		
		\begin{equation}
		d(X,Y):=\KL(X\|Y)+\KL(Y\|X).
		\end{equation}
	\end{definition}
\begin{remark}
			Since the KL-divergence can be defined by  distribution functions (or random vectors), so does $d$. Specifically, given two distribution $P, Q$,  $d(P,Q):=\KL(P\|Q)+\KL(Q\|P)$. Clearly, if $X, Y$ has distribution $P_X, P_Y$, then $d(X,Y)=d(P_X,P_Y)$. 
		\end{remark}
	
	To the best of our knowledge, \cite{SKL-spatiogram} showed the first application of the symmetric KL-divergence as a spatiogram similarity measure  in 2011, whose idea was motivated by Vasconcelos (2004)\cite{NV}, and then the symmetric KL-divergence has been used in plenty of fields until now.
	As discussed by Domke (2021)\cite{SKL-diagnostic}, when it comes to interpreting diagnostic inference problems, certain methods like Laplace's method prove impractical without the inclusion of sample distribution information, in addition to a mere set of samples. Moreover, the KL-divergence measure provides limited insight into the absolute error present in approximate algorithms. However, the symmetric KL-divergence overcomes these shortcomings, thus serving as a valuable error measure in simulation-based interpret diagnostics. Andriamanalimanana (2019)\cite{SKL-anomaly} concluded that   the inherent symmetry advantages of the symmetric KL-divergence, when compared to the KL-divergence, further extends its application as a numerical score within anomaly detection systems. In the field of text classification, Chen \textit{et al.} (2018)\cite{SKL-Centroid}  discovered that centroid estimation based on the symmetric KL-divergence outperforms the conventional naive Bayes method, which employs average estimated centroids or relies solely on KL-divergence.

	In the realm of information theory, symmetric KL-divergence has a strong connection with channel theory. Aminian \textit{et al}.(2014)\cite{SKL-Capacity1} demonstrated that for any channel $W=P(y|x)$, where $x$ represents the input signal with alphabet $\mathcal{X}$,and $y$ represents the output with alphabet $\mathcal{Y}$, the capacity $C$ of $W$ is bounded above by the expression:
		$$C\le \max_{P_x\in\mathcal{P}(\mathcal{X})} d(P_{xy}, P_{x}\times P_y),$$
		where $\mathcal{P}(\mathcal{X})$ is the set of all distribution on $\mathcal{X}$, $P_{x}$ means a distribution of  the input signal and $P_{y}$ means a distribution of  the output signal, $P_{xy}$ means the joint distribution of the signals (Since channel $W$ is fixed, $P_y$ is uniquely determined by $P_x$). 	
		Specifically, in the case of an add white Gaussian noise channel with a fixed signal to noise ratio (SNR), we have
		\begin{align}
		C&=\frac{1}{2}\log(1+\mbox{SNR})\le \max_{P_x\in\mathcal{P}(\mathcal{X})} d(P_{xy}, P_{x}\times P_y)\\
		&=d(P_{xy}, P_{x}\times P_y)=\V(X_{in})/\sigma^2=\mbox{SNR},\;\;\forall P_x \in \mathcal{P}(\mathcal{X}),
		\end{align}
		where $X_{in}$ is the random variable with distribution $P_x$,
	$\sigma^2$ is the variance of the channel noise. These observations are utilized in diffusion-based molecular communication for capacity estimation, as discussed in \cite{SKL-Capacity1} and \cite{SKL-Capacity2}. Furthermore, the symmetric KL-divergence is also called as the bidirectional KL-divergence in neural networks and deep learning fields,   considered as a loss function in the corresponding optimization model recently (\cite{R-drop}, \cite{ruan2023gist}).

	The combination of Gaussian distribution and symmetric KL-divergence has also been explored in various contexts, such as \cite{SKL-Wasserstein} and \cite{SKL-spatiogram}. For instance, Welandawe \textit{et al}.(2022)\cite{SKL-Wasserstein} employed symmetric KL-divergence to detect convergence in Black-box variational inference in Machine Learning. They demonstrated that the Wasserstein distance, often used in MCMC diagnostic for detecting stationarity, can be controlled by the symmetric KL-divergence under the exponentially controlled assumption. Given that this assumption holds when discussing the symmetric KL-divergence $d(X,\hat{G})$, where $\hat{G}$ is Gaussian, Welandawe \textit{et al}. focused on the symmetric KL-divergence for two different Gaussian random variables and deduced the convergence order as the two distributions gradually approach each other.

	In summary, symmetric KL-divergence is a common measure that indicates the similarity between two distributions. Therefore, it is natural to consider symmetric KL-divergence as a standard for measuring the distance of a random variable with Gaussian distribution.  According to the definition and Pinsker's inequality \eqref{Pinsker}, we have:
	\begin{equation}\label{SKL is stronger}
	d(X, G)\ge \KL(X\|G),\;\; d(X,G)\ge \|P_X-\Phi\|^2_{TV},
	\end{equation}
	where $G$ denotes the standard Gaussian random variable, $\Phi$ is the distribution of $G$ and the notations will be used throughout the paper.	Consequently, the CLT under symmetric KL-divergence   implies a much stronger convergence compared to KL-divergence, total variation distance, or Wasserstein distance, as stated in     \cite{SKL-Wasserstein}.

	By discussing the CLT under symmetric KL-divergence, we not only establish a uniform upper bound for the aforementioned convergence measures, but also cater to situations that specifically require symmetric KL-divergence as a similarity metric, since its wide range of applications. Despite some research being conducted on the convergence properties of symmetric KL-divergence, such as \cite{SKL-spatiogram} and \cite{SKL-Wasserstein}, little attention has been given to the CLT or other probability convergence theorems. This could be attributed to the challenging analysis of the term $KL(G\|X)$. Specifically, this term cannot be easily expressed in terms of entropy like $KL(X\|G_X)=h(G_X)-h(X)$, and requires X to have a strictly positive density $p(x)>0$, for all $x\in \R$.

	In this paper, we use Stein's method (as introduced in Section \ref{section2}) to overcome the  untractable problem of symmetric KL-divergence. Since we aim to deduce the CLT under symmetric KL-divergence, the paper only focus on
	\begin{equation}
	d(W_n, G)=\KL(W_n\|G)+\KL(G\|W_n),
	\end{equation}
	where $W_n$ is a normalized sum of independent random variables with $\E W_n=0$ and  $\V(W_n)=1$.  Given that following regular assumption  is valid:
	\begin{assumption}\label{part A assume}
		Given  random variables  $X_1, X_2,\cdots, X_n$ with density functions $\rho_1, \rho_2, \cdots, \rho_n$. Assume that
		$\exists \delta_1\in(\frac{1}{4},1)$ and set $A_n\subset\{1,2,\cdots, n\}$ \textit{s.t} $\delta_1\le\frac{|A_n|}{n}$ and
		\begin{equation}
		\rho_i(x)\ge l_1e^{-l_2\frac{x^2}{2}},\forall x\in \R, \forall i\in A_n,
		\end{equation}
		for some $l_1>0, l_2>0$.
	\end{assumption}
	Some remarks on Assumption \ref{part A assume} are discussed in Section \ref{section2}.  Our main result can be summarized as follows.
	\begin{theorem}\label{bound for d under Stein}
		Given a random variable sequence  $X_1, X_2,\cdots, X_n, \cdots$ \textit{s.t.} for any $i=1,2,\cdots,$ $\E X_i=0, \sup_i J(X_i)<J<\infty, \sup_i \E|X_i|^{4+\delta_0}<M<\infty$ for some $\delta_0>0$.
		Let
		$W_n$ be a normalized sum of independent random variables, \textit{i.e.}
		$$W_n=\frac{X_1+X_2+\cdots +X_n}{\sqrt{\sum_{i=1}^n\V(X_i)}}.$$
		If Assumption \ref{part A assume} holds for $X_1, X_2,\cdots, X_n$ with $(\delta_1, l_1,l_2)$ when $n$ is large enough,
		%		and
		%		\begin{equation}\label{??}
		%		 k>\frac{2}{\log_M(4\delta_1)}-2-\delta_0,
		%		\end{equation}
		%			\begin{equation}\label{??2}
		%		k\left( 1-\frac{M^{\frac{2}{4+delta_0}}}{4\delta_1}\right) >\frac{15}{8},
		%		\end{equation}
		then there exists a
		constant $C=C(J,M,\delta_0,\delta_1,l_1,l_2)$ such that the symmetric KL-divergence is bounded
		by
		\begin{equation}
		d(W_n, G)\le C\frac{\ln n}{\sqrt{n}},
		\end{equation}
		where $G$ is the standard Gaussian random variable.
	\end{theorem}

	Additionally, the well-known local central limit theorem says that the total variation distance between $W_n$ and standard Gaussian random variable $G$ meets the order $O(1/\sqrt{n})$ (\textit{e.g.} Lemma \ref{p-phi to 1/n^k}). Thanks to  \eqref{SKL is stronger}, we conclude that the convergence rate of $d(W_n, G)$ is at least  $O(1/n)$.
	However,  Theorem \ref{bound for d under Stein} still   remains a gap of the rate estimation between $O(\ln n/\sqrt{n})$ and $O(1/n)$. In order to overcome this gap, we add the finite exponential moment condition and combine the result in \cite{Esseenbound}, and then we obtain the 'best rate estimation' $O(1/n)$ for symmetrical KL-divergence $d(W_n, G)$ (see Section \ref{section4} for details).
	
	\textbf{Paper outline:} In Section \ref{section2}, we introduce the Stein's method and construct our primary model for symmetrical KL-divergence in the form of the Stein equation. Meanwhile, a core Property \ref{e^e assumption}  is discussed, which is a more general condition of Theorem \ref{bound for d under Stein} and  replace the Assumption \ref{part A assume} in our proof.  The local limit theorem, which plays a crucial role in our findings, is   presented in the final of  Section \ref{section2}. The proof and additional observations regarding Theorem \ref{bound for d under Stein} are presented in Section \ref{section3}. In Section \ref{section4}, some assumptions to further improve the rate in  Theorem \ref{bound for d under Stein} are discussed.
	
	\section{Basic Ideas}\label{section2}

	\subsection{Stein's Method}
	Stein (1972)\cite{Stein1972ABF} discovered   that a random variable $\hat{G}$ follows a Gaussian distribution $\hat{G}\sim\mathcal{N}(0,\sigma^2)$ if and only if
	\begin{equation}\label{stein idea}
	\E(\hat{G} f(\hat{G}))=\sigma^2\E(f'(\hat{G}))
	\end{equation}
	for any absolutely continuous function $f$ for which the above expectation exists. Hence, \eqref{stein idea} can be regarded as a criterion for the "Gaussianity" of a random variable. Specifically, when $|\E(X f(X))-\sigma^2\E(f'(X))|$ is small for any absolutely continuous function $f$, we may infer that $X$ is close to   the Gaussian random variable $\hat{G}$.
	
	Without loss of generality, standard Gaussian random variable $G$ is considered after now. Building upon \eqref{stein idea} and Taylor expansion,   Goldstein and   Reinert (1997) \cite{Goldstein1997SteinsMA} presented the ensuing inequality for any absolutely continuous function $g$ with $\|g'\|_\infty<\infty$ and a random variable $W$ with $\E W=0$ and $\V(W)=1$:
	\begin{align}
	|\E g(W)-\E g(G)|&=|\E(f_g'(W)-Wf_g(W))|\\
	&=|\E(f_g'(W)-f'_g(W^*))|\\
	&\le \|f_g''\|_\infty \E|\Delta|,\label{delta}
	\end{align}
	where $W^*$ represents the zero bias\footnote{Given a random variable $V$ with mean zero and variance $\sigma^2$, we say $V^*$ is  the zero bias of $V$ if $\sigma^2\E f'(V^*)=\E[Vf(V)]$ for all absolutely continuous function $f$ for which these expectations exist.} of $W$,   $\Delta:=W^*-W$, and $f_g$ is the solution of Stein's equation with a given function $g$ as follows:
		\begin{equation}\label{stein equation}
		f_g'(w)-wf_g(w)=g(w)-\E g(G).
		\end{equation}
	It should be noted, as proved by Stein (1986) \cite[Lemma 3]{stein1986} (see also Rai$\check{\mbox{c}}$ (2004) \cite[Lemma 1]{Rai2004AMC}), that upper bounds exist for the solution $f_g$ of Stein's equation generated by $g$, as follows.
	\begin{lemma}\label{lemma of basic bound in stein}
		If $g$ is an absolutely continuous function, then
		\begin{equation}\label{basic bound in stein}
		\|f_g\|_\infty\le 2\|g'\|_\infty,\;\; \|f_g'\|_\infty\le \sqrt{2/\pi}\|g'\|_\infty,\;\; \|f_g''\|_\infty\le 2\|g'\|_\infty.
		\end{equation}
		where $f_g$ is the solution of Stein's equation with respect to the function $g$, as \eqref{stein equation}.
	\end{lemma}
	Therefore,  we can promptly obtain a CLT using Stein's method, which can be summarized as follows.
	\begin{lemma}
		For any absolutely continuous function $g$ with $\|g'\|_\infty<\infty$ and random variable $W$ with $\E W=0$ and $\V(W)=1$:
		\begin{align}
		|\E g(W)-\E g(G)|\le 2\|g'\|_\infty \E|\Delta|,
		\end{align}
		where $G$ denotes the standard Gaussian random variable, $W^*$ represents the zero bias of $W$, and $\Delta:=W^*-W$.
		In particular, when $W=\sum_{i=1}^n \xi_i$ represents the sum of independent random variables $\xi_1, \xi_2,\cdots \xi_n$ with $\E \xi_i=0$ for all $i$ and $\sum_{i=1}^n\V(\xi_i)=1 $, we have
		\begin{align}
		|\E g(W)-\E g(G)|\le 3\|g'\|_\infty \sum_{i=1}^n \E|\xi_i|^3.
		\end{align}
	\end{lemma}
	The proof can be found in \cite[Section 3]{Goldstein1997SteinsMA}(see also \cite[Corollary 1]{stein1986} and \cite[Theorem 3.1]{steinbook}).
	
	Thanks to the feature abstraction idea, as \eqref{stein idea} for Gaussian random variable, Stein’s method can not only address asymptotic Gaussianity for independent sum cases, but it can also be applied to more general cases and even to cases with other limiting distributions, such as the Poisson distribution, as \eqref{stein idea} becomes  other appropriate expressions.  Despite being a powerful tool, Stein's method still encounters challenges when the function $g$ lacks a bounded first-order derivative, and $g$ may depend on the step parameter $n$, which tends to infinity.
	
	In summary, Stein's method extracts the feature of a Gaussian random variable that \eqref{stein idea} holds, while the entropy method may note that a Gaussian random variable has maximum entropy when its variance is fixed. These two ideas  capture the characteristics of the Gaussian distribution and then deduce a CLT, thus it is natural to consider them together.

	\subsection{Symmetrical KL-divergence in Stein Form}

	In order to ensure $\KL(G\|X)$ in symmetric KL-divergence be well-defined,   we always consider a random variable with strictly positive and absolutely continuous density unless otherwise stated.
	
	Given a random variable $X$ with mean $0$, variance $1$ and density $p$, denote $g=\ln p$, then we directly obtain
	\begin{equation}\label{model}
	\begin{aligned}
	d(X,G)&=\int_\R [p(x)\ln p(x)-p(x)\ln \phi(x) +\phi(x)\ln\phi(x) -\phi(x)\ln p(x)]  dx\\
	&=\E g(X)-\E g(G),
	\end{aligned}
	\end{equation}	
	where $\phi$ is the density function of standard Gaussian random variable $G$. Replace $X$ by $W_n$ defined in Theorem \ref{bound for d under Stein}, we have
	\begin{equation}
	g'(x)=\rho_{W_n}(x)=\frac{p'_n(x)}{p_n(x)},
	\end{equation}
	where $p_n$ is the density of $W_n$, and $\rho_{W_n}$ is the score function of $W_n$. Undoubtedly, it is impossible to impose a uniform bound on the score functions $\rho_{W_n}(x)$, hence we cannot  apply Lemma \ref{lemma of basic bound in stein} directly. However, the local limit theorem presented below can address the issue of unbounded derivatives in our Stein function model \eqref{model}.  Note that for the rest of this paper, we will abbreviate $d(W_n,G)$ as $d(W_n)$.
	
	\subsection{A Core Property of the Density of $W_n$}
	Clearly, to deal with the convergence of $d(W_n)$, we need to focus on the density function $p_n$ of the partial sum $W_n$. For the proof of Theorem \ref{bound for d under Stein}, we weaken Assumption \ref{part A assume} by the following property of the density of partial sum.
	\begin{property}\label{e^e assumption}
		Given the density $p_n$ of $W_n$, where
		$$W_n=\frac{X_1+X_2+\cdots +X_n}{\sqrt{\sum_{i=1}^n\V(X_i)}},$$
		%$\E X_i=0, \sup_i J(X_i)<J<\infty, \sup_i \E|X_i|^{4+\delta_0}<M<\infty$ for some $\delta_0>0, \forall i=1,2,\cdots$.
		there exists $a>1, k_1>0, k_2>0, s>0,$ and $v\ge 0$, s.t.
		\begin{equation}\label{k>v}
		\frac{a}{a-1}\Big( \frac{1}{2}+s\Big) <2,
		\end{equation}
		and
		\begin{equation}\label{p>e}
		p_n(x)\ge k_1e^{-n^sk_2e^{\frac{x^2}{2a} }},\;\;\forall |x|>v\sqrt{\ln n}.
		\end{equation}
	\end{property}
	\begin{remark}
		Obviously,
		Property \ref{e^e assumption} can be written as limit form: for the density $p_n$ of $W_n$, we assume that $\exists a>1, s>0$, s.t. for any large $n$,
		\begin{equation}
		\liminf_{|x|\to\infty}  \frac{e^{-\frac{x^2}{2a}} \ln p_n(x)}{n^s}>-\infty.
		\end{equation}
	\end{remark}
	
	We have proven that under the Assumption \ref{part A assume}, Property \ref{e^e assumption} holds for $W_n$ immediately (see  \ref{App1}),
		henceforth, 	we only focus on the Property \ref{e^e assumption} rather than Assumption \ref{part A assume}. A remark for Property \ref{e^e assumption} is discussed as follow.
	
	\begin{remark}
		In order to make
		$d(W_n)<\infty,$ we need the term
		$\int_\R \phi(x) \ln p_n(x)dx$ to be finite, \textit{i.e.}
		\begin{equation}\label{L1 integrable for ln p}
		\E |\ln p_n(G)|<\infty.
		\end{equation}
		Since $|\phi \ln p_n|\le p_n\phi\le p_n\in L^1(\R)$ when $\ln p_n>0$, we only need a lower bound for $p_n$. Note that if on an infinite Lebesgue measure set $S$,
		$$p_n(x)\le e^{-e^{x^2/2}},\;\;\forall x\in S,$$
		then $\E |\ln p_n(G)|\ge \int_S e^{x^2/2}\phi(x)dx=\infty.$
		Thus  to a certain extent, Property \ref{e^e assumption}  is considered as a quiet  necessary condition  for \eqref{L1 integrable for ln p}.
	\end{remark}
	
	\subsection{Local Limit Theorem}
	The local limit theorem elucidates the distance between the density function $p_n$ of $W_n$ and the density function $\phi$ of a standard Gaussian random variable $G$. By virtue of its meticulous expansion,   $p_n$ can be approximated by $\phi$ and thereby attain a more manageable scenario. In this section, we present the local limit theorem that is utilized in our proof for Theorem \ref{bound for d under Stein}.

	Using the classical results in \cite{Sumofindependent}: Theorem 14 in Chapter VII, (1.10) in Chapter VI, and the statement in Page 173, and (6.13) in \cite{1986Normalapproximation}, we conclude that,
	\begin{lemma}\label{p-phi to 1/n^k}
		Consider a sequence  of independent random variables $X_1,X_2,\cdots,$ with densities $\rho_{1},\rho_{2},\cdots,$ and $\E X_i=0, \forall i=1,2,\cdots$. Denote
		\begin{equation}
		B_n:=\sum_{i=1}^n \V(X_i).
		\end{equation}
		If the following conditions hold:
		\begin{enumerate}
			\item $\liminf_{n\to\infty}\frac{B_n}{n}>0.$
			
			\item $\exists k\in\mathbb{N}^+, \delta_0>0, \limsup_{n\to\infty}\frac{1}{n}\sum_{i=1}^n \E |X_i|^{k+2+\delta_0}<\infty.$
			
			\item There exists a subsequence $X_{n_1}, X_{n_2}\cdots, X_{n_m},\cdots$, s.t. for some $\lambda>0$,
			\begin{equation}
			\liminf_{n\to\infty} \frac{\#\{m:n_m\le n  \}}{n^\lambda}>0,
			\end{equation}
			and $\sup_{n_m} \|\rho_{{n_m}}\|_{TV}<\infty$, where $\#$ represents the counting measure, and $ \|\rho\|_{TV}$ is the total variation of function $\rho$ defined as
				\begin{equation*}
				\|\rho\|_{TV}:=\sup_{[a,b]\subset \R}\sup_{P\in\mathcal{P}_{ab}}\sum_{i=1}^{n_P-1}|\rho(x_i)-\rho(x_{i-1})|,
				\end{equation*}
				$\mathcal{P}_{ab}=\{P=\{x_0,x_2,\cdots,x_P\}|P\mbox{~is a partition of~}\; [a,b]\}$.
		\end{enumerate}
		Then the density $p_n(x)$ of the random variable $W_n:=B_n^{-1/2}\sum_{i=1}^n X_i$ satisfies that
		\begin{equation}
		p_n(x)-\phi(x)=\sum_{i=1}^k \frac{1}{n^{i/2}}R_i(x)+o\left(\frac{1}{n}\right),
		\end{equation}
		where $\phi(x)=(\sqrt{2\pi})^{-1}e^{-x^2/2}$ is the density of standard normal distribution, all $R_i(x)$ have the form $R_i(x)=H_i(x)\phi(x)$, $H_i(x)$s are polynomial,
		and $o\left(\frac{1}{n}\right)$ is uniformly about $x$.
	\end{lemma}
	
	\begin{remark}
		Specifically, for example
		\begin{equation}
		R_1(x)=\left(\frac{1}{6}\frac{\sqrt{n}}{B_n^{1/2}}\frac{1}{B_n}\sum_{i=1}^n \gamma_{3,i}\right)(x^3-3x)\phi(x),
		\end{equation}
		\begin{equation}
		R_2(x)=\frac{1}{72}\left(\frac{\sqrt{n}}{B_n^{3/2}}\sum_{i=1}^n \gamma_{3,i}\right)^2(x^4-6x^2+3)\phi(x)+\left(\frac{n}{24B^2_n}\sum_{i=1}^n \gamma_{4,i}\right)(x^3-3x)\phi(x),
		\end{equation}
		where
		\begin{equation}
		\gamma_{3,i}=\E X_i^3,
		\end{equation}
		and
		\begin{equation}
		\gamma_{4,i}=\E X_i^4-3\E X_i^2.
		\end{equation}
	\end{remark}
	
	Note that 	for a random variable $Y$ with strictly positive and absolutely continuous density function $\rho$, we always have
	\begin{equation}\label{TV and c.h.f}
	\V(Y)J(Y)\ge 1 ,\;\; \|\rho\|_{TV}\le \int_\R |\rho'|=\int_\R \sqrt{\rho}\frac{|\rho'|}{\sqrt{\rho}}  \le 1\cdot \sqrt{\int_\R \frac{|\rho'|^2}{\rho}}=\sqrt{J(Y)},
	\end{equation}
	\textit{i.e.} finite Fisher information guarantees the  non-degenerate variance  condition 1 and the bounded total variation condition 3  hold simultaneously in Lemma \ref{p-phi to 1/n^k}. Based on Lemma \ref{p-phi to 1/n^k},  we turn to the proof of our main result Theorem \ref{bound for d under Stein}.

	\section{The Proof of Theorem \ref{bound for d under Stein}}\label{section3}
	This section we show the proof of our main result Theorem \ref{bound for d under Stein} with three steps.
	
	\subsection{Applying Stein’s method}
	
	We use a more precise bound for the Stein equation to deal with our unbounded derivative case. Specifically,
	according to the Appendix of Chapter 2 in \cite{steinbook}, we can give a bound with the third order Taylor expansion for {\it even function} instead of Lemma \ref{lemma of basic bound in stein}. For $\Delta$   defined in \eqref{delta} and let $\xi$ be a random variable between $0$ and $1$,   then by third order
		Taylor expansion one gets
	\begin{align}
	|\E h(W_n)- \E h(G)|&=|\E\left( f_h'(W_n)-f_h'(W_n+\Delta)    \right)|\\
	&=\left|\E\left( (f_h''(W_n)+\frac{1}{2}f_h'''(W_n+\xi\Delta)\Delta)\Delta \right)\right|\label{3.2}\\
	&\le \sqrt{\E(f_h''(W_n))^2}\sqrt{\E\Delta^2}+\|f_h'''\|_\infty\E\Delta^2, \label{term E D^2}
	\end{align}
	where $h$ is an even function. In order to apply Stein’s method successfully, we  need to bound $|f''_h|$ and $|f_h'''|$, which is accomplished as follows.
	\begin{equation}\label{bound for fh''}
	\begin{aligned}
	&|f_h''(x)|
	\le|h'(x)|+|(\sqrt{2\pi}(1+x^2)e^{x^2/2}(1-\Phi(x))-x)\int_{-\infty}^x  h'(t)\Phi(t) dt\\
	&\;\;\;\;+(\sqrt{2\pi}(1+x^2)e^{x^2/2}\Phi(x)+x)\int_x^\infty   h'(t)(1-\Phi(t)) dt|\\
	&=|h'(x)|\\
	&\;\;\;\;+\bigg{|}(\sqrt{2\pi}(1+x^2)e^{x^2/2}(1-\Phi(x))-x)\left( h(x)\Phi(x)-\int_{-\infty}^x  h(t)\phi(t) dt \right)\\
	&\;\;\;\;+(\sqrt{2\pi}(1+x^2)e^{x^2/2}\Phi(x)+x)\left( \int_{x}^\infty  h(t)\phi(t) dt-h(x)(1-\Phi(x)) \right)\bigg{|}\\
	&\le|h'(x)|+2\sqrt{2\pi}(1+x^2)e^{x^2/2}\Phi(-|x|)\left(\E h(G)+ \frac{1}{\Phi(-|x|)}\int_{-\infty}^{-|x|} h(t)\phi(t)dt\right)\\
	&\;\;\;\;+|x\E h(G)|+|xh(x)|\\
	%&\le|h'(x)|+(2\sqrt{2\pi}(1+x^2)e^{x^2/2}\Phi(x)(1-\Phi(x))+x)\|h\|_\infty+|xh(x)|\\
	&\le |h'(x)|+C_0|x\E h(G)|+C_0\bigg|x(\frac{1}{\Phi(-|x|)}\int_{-\infty}^{-|x|} h(t)\phi(t)dt)\bigg|+|xh(x)|,
	\end{aligned}
	\end{equation}
	and
	\begin{align}
	|f_h'''(x)|&=|2f_h'(x)+xf_h''(x)+h''(x)|\\
	&\le 2(\sqrt{2/\pi}+|x|)\|h'\|_\infty+\|h''\|_\infty.\label{bound for fh'''}
	\end{align}

		We also need to estimate the term $\E\Delta^2$ in \eqref{term E D^2}.
		\begin{lemma}\label{bound of delta^2}
			Under the the conditions of Theorem \ref{bound for d under Stein}, we conclude that
			\begin{equation}
			\E\Delta^2\le \gamma \frac{1}{n},
			\end{equation}
			where 
			\begin{equation}\label{gamma}
			\gamma=\frac{4}{3}J^2M^{\frac{4}{4+\delta_0}}.
			\end{equation}
		\end{lemma}
		\pf
		According to \cite[Lemma 2.8]{steinbook}, we have  
		\begin{equation}
		\Delta=\xi_I^*-\xi_I,
		\end{equation}
		where $\xi_i=\frac{X_i}{\sqrt{\sum_{j=1}^n \V(X_i)}}$, $\xi_i^*$ is a zero bias of $\xi_i$ and $I$ is a random index which is  independent of $\xi_i, \xi_i^*, i=1,2,\cdots,n$ and has the following distribution
		\begin{equation}
		\P(I=i)=\V(\xi_i)=\frac{\E X_i^2}{\sum_{j=1}^n \V(X_i)}.
		\end{equation}
		Thus
		\begin{align}
		\E\Delta^2&=\sum_{i=1}^n \frac{\E X_i^2}{\sum_{j=1}^n \V(X_i)}\E|\xi_i^*-\xi_i|^2\\
		&\le 2\sum_{i=1}^n \frac{\E X_i^2}{\sum_{j=1}^n \V(X_i)}(\E(\xi_i^*)^2+\E(\xi_i^2))\\
		&\overset{(a)}{=}2\sum_{i=1}^n\left( \frac{1}{3}\E \xi_i^4+(\E \xi_i^2)^2\right) \\
		&\overset{(b)}{\le} \frac{4}{3} \sum_{i=1}^n\E \xi_i^4\\
		&=\frac{4}{3(\sum_{j=1}^n \V(X_i))^2}\sum_{i=1}^n\E X_i^4\\
		&\overset{(c)}{\le} \frac{4J^2}{3n^2} \sum_{i=1}^n(\E X_i^{4+\delta_0})^{\frac{4}{4+\delta_0}}\\
		&\le \frac{4}{3n}J^2M^{\frac{4}{4+\delta_0}},
		\end{align}
		where $(a)$ comes from the definition of zero bias, $(b)$ holds by Jensen's inequality, and $(c)$  holds according to   CR bound  $\V(\xi)J(\xi)\ge 1$ for any random variable and Jensen's inequality. 
		\e

		Therefore, we conclude that
		\begin{equation}
		|\E h(W_n)- \E h(G)|=\sqrt{\gamma\E(f_h''(W_n))^2}\sqrt{\frac{1}{n}}+\gamma\|f_h'''\|_\infty\frac{1}{n},\label{Stein f'''}
		\end{equation}
		where $\gamma$ is independent of $n$.
	
	These bounds will be used in the proof.
	
	\subsection{Truncation}\label{section 3.2}
	Since we assume  $\sup_i \E |X_i|^{4+\delta_0}<M<\infty, \sup_i J(X_i)<J<\infty $,  according to  Lemma \ref{p-phi to 1/n^k} (with $k=2$) we  conclude that for any sufficient large $n$,
	\begin{equation}\label{local bounds}
	-\frac{C}{n}-C\frac{x^4+1}{\sqrt{n}}\phi(x)\le p_n(x)-\phi(x)\le \frac{C}{n}+C\frac{x^4+1}{\sqrt{n}}\phi(x),\;\;\forall x\in\R,
	\end{equation}
	where $C$ is a constant that only depends on $M$. We denote $r(x):=\frac{C}{n}+C\frac{x^4+1}{\sqrt{n}}\phi(x)$ ($\le 1-\phi(x)$ when $n$ is large enough).

	Note that on $A=\{x; p_n(x)\ge \phi(x)\},$ $\ln p_n\le \ln( \phi(x)+r(x))$;
	on $A^c=\{x; p_n(x)< \phi(x)\},$ $\ln p_n\ge  \ln( \phi(x)-r(x))$, and
	$\ln( \phi(x)+r(x) )(p_n(x)-\phi(x))\ge 0.$
	Thus by denoting $B(u):=\{ x; |x|\le u\sqrt{\ln n}\},$ $0<u<\sqrt{2},$   when $n$ is large we have
	\begin{equation}\label{p-r}
	\phi(x)-r(x)\ge \frac{1}{\sqrt{2\pi}}\frac{1}{n^{u^2/2}}\Big(1-C\frac{x^4+1}{\sqrt{n}}\Big)-\frac{C}{n}\ge \frac{1}{n},\;\;\forall x\in B(u), 0<u<\sqrt{2},
	\end{equation}
	and
	\begin{align*}
	&d(W_n)=\int_{B(u)}(p_n(x)-\phi(x)) \ln p_n(x) dx +\int_{B(u)^c}(p_n(x)-\phi(x)) \ln p_n(x) dx\\
	&=\int_{A\cap B(u)}(p_n(x)-\phi(x)) \ln p_n(x) dx+\int_{A^c\cap B(u)}(p_n(x)-\phi(x)) \ln p_n(x) dx\\
	&\;\;\;\; +\int_{B(u)^c}(p_n(x)-\phi(x)) \ln p_n(x) dx\\
	&\le \int_{A\cap B(u)}(p_n(x)-\phi(x))\ln( \phi(x)+r(x))dx+\int_{A^c\cap B(u)}(p_n(x)-\phi(x)) \ln( \phi(x)-r(x))dx\\
	&\;\;\;\;+\int_{B(u)^c}(p_n(x)-\phi(x)) \ln p_n(x) dx\\
	&\le \left( \int_{A\cap B(u)}(p_n(x)-\phi(x))\ln\left(  \frac{\phi(x)+r(x)}{\phi(x)-r(x)}\right)  dx\right)\\
	&\;\;\;\;+\left( \int_{A\cap B(u)}(p_n(x)-\phi(x))\ln( \phi(x)-r(x))dx+ \int_{A^c\cap B(u)}(p_n(x)-\phi(x))\ln( \phi(x)-r(x))dx \right)\\
	&\;\;\;\;+\int_{B(u)^c}(p_n(x)-\phi(x)) \ln p_n(x) dx\\
	&\le I_1+I_2+I_3+I_4,
	\end{align*}
	where
	$$I_1=\int_{A\cap B(u)}(p_n(x)-\phi(x))\ln\left(  \frac{\phi(x)+r(x)}{\phi(x)-r(x)}\right) dx ,$$
	$$I_2=\E h_1(W_n)-\E h_1(G),$$
	$$I_3=\int_{B(u)^c}|(p_n(x)-\phi(x))h_1(x)|dx,$$
	$$I_4=\int_{B(u)^c}(p_n(x)-\phi(x)) \ln p_n(x) dx,$$
	\begin{equation}
	h_1=\begin{cases}
	[(\ln n)^2(x-u\sqrt{\ln n})+\ln(\phi(u\sqrt{\ln n})-r(u\sqrt{\ln n}))]^-, &x>u\sqrt{\ln n},\\
	\ln( \phi(x)-r(x)), & |x|\le u\sqrt{\ln n},\\
	[-(\ln n)^2(x+u\sqrt{\ln n})+\ln(\phi(-u\sqrt{\ln n})-r(-u\sqrt{\ln n}))]^-, &x<-u\sqrt{\ln n},
	\end{cases}
	\end{equation}
	and $[x]^-:=x\wedge 0$. Clearly, $h_1$ is an even function.

	\begin{remark}
		We also need some smoothing for $h_1$ on four little sets $(-u\sqrt{\ln n}-\varepsilon, -u\sqrt{\ln n}), (u\sqrt{\ln n}, u\sqrt{\ln n}+\varepsilon)$ and the neighborhoods of two zero arriving points in outer area $|x|>u\sqrt{\ln n}$ . But this do not change any following bounds, so we just write $h_1$ with the above form.
	\end{remark}
	
	\subsection{The Remaining Proof}
	
	Fix $1\le u<\sqrt{2}$ ($u=\sqrt{1.5}$ for example), We estimate $I_i, i=1,2,3,4$ one by one.
	\begin{enumerate}
		\item Upper bound of $I_1$.  For $n$ large enough,
		\begin{align}
		I_1&=\int_{A\cap B(u)}(p_n(x)-\phi(x))\ln\left(  \frac{\phi(x)+r(x)}{\phi(x)-r(x)}\right) dx\\
		&=\int_{A\cap B(u)}(p_n(x)-\phi(x))\ln\left(  \frac{1+\frac{C(x^4+1)}{\sqrt{n}}+\frac{C}{n}\phi(x)^{-1}}{1-\frac{C(x^4+1)}{\sqrt{n}}-\frac{C}{n}\phi(x)^{-1}}\right) dx\\
		&\le \int_{A\cap B(u)}(p_n(x)-\phi(x))\ln\left(  \frac{1+\frac{C(x^4+1)}{\sqrt{n}}+\frac{C}{n}\sqrt{2\pi}n^{u^2/2}}{1-\frac{C(x^4+1)}{\sqrt{n}}-\frac{C}{n}\sqrt{2\pi}n^{u^2/2}}\right) dx\\
		&\le \int_{A\cap B(u)}(p_n(x)-\phi(x)) \left( \frac{1+\frac{C(4(\ln n)^{2}+1)}{\sqrt{n}}+\frac{C}{n}\sqrt{2\pi}n^{u^2/2}}{1-\frac{C(4(\ln n)^{2}+1)}{\sqrt{n}}-\frac{C}{n}\sqrt{2\pi}n^{u^2/2}}-1\right)dx\\
		&\le \int_{A}(p_n(x)-\phi(x)) \left( \frac{1+\frac{C(4(\ln n)^{2}+1)}{\sqrt{n}}+\frac{C}{n}\sqrt{2\pi}n^{u^2/2}}{1-\frac{C(4(\ln n)^{2}+1)}{\sqrt{n}}-\frac{C}{n}\sqrt{2\pi}n^{u^2/2}}-1\right)dx\\
		&\le C_1\left( \frac{1}{n^{(2-u^2)/2}}+\frac{(\ln n)^{2}}{\sqrt{n}}\right) \|p_n-\phi\|_1.
		\end{align}
		 Note that
		\begin{equation}
		\|p_n-\phi\|_1=2\int_{A^c} \phi(x) -p_n(x) dx\le 2\int_{|x|\le \ln n} |\phi(x) -p_n(x)|dx+2\int_{|x|>\ln n}\phi(x)dx,
		\end{equation}
		according to Lemma \ref{p-phi to 1/n^k}, it is easy to check that under the conditions of Theorem \ref{bound for d under Stein}, $\|p_n-\phi\|_1=O(\frac{\ln n}{\sqrt{n}})$,
		thus we have
		\begin{equation}\label{I_1,2}
		I_1\le C_2\frac{1}{\sqrt{n}}.
		\end{equation}

		\item Upper bound of $I_2$.  Note that for any $x\in B(u)$ and large $n$,
		\begin{equation}\label{1+C}
		\begin{aligned}
		&|(\ln(\phi(x)-r(x))'|\le|x|\frac{(1-C(x^4+1)/\sqrt{n})\phi(x)}{\phi(x)-r(x)}+\frac{8x^{3}C/\sqrt{n}\phi(x)}{\phi(x)-r(x)}\\
		&=|x|\frac{(1-C(x^4+1)/\sqrt{n})\phi(x)-\frac{C}{n}+\frac{C}{n}}{(1-C(x^4+1)/\sqrt{n})\phi(x)-\frac{C}{n}}+\frac{8x^{3}C/\sqrt{n}\phi(x)}{\phi(x)-r(x)}\\
		&=|x|+|x|\frac{\frac{C}{n}}{\phi(x)-r(x)}+\frac{8x^{3}C/\sqrt{n}\phi(x)}{\phi(x)-r(x)}\\
		&\overset{(a)}{\le}|x|+|x|\frac{\frac{C}{n}}{\frac{1}{n}}+\frac{8x^{3}C/\sqrt{n}\phi(x)}{\phi(x)-r(x)}\\
		&\overset{(b)}{\le} (1+C)|x|+\frac{(1-C(x^4+1)/\sqrt{n})\phi(x)}{\phi(x)-r(x)}\\
		&\le (1+C)(|x|+1),
		\end{aligned}
		\end{equation}
		where we use \eqref{p-r} to deduce $(a)$, and $(b)$ is from the fact that
		\begin{equation}
		\frac{8x^{3}C }{\sqrt{n}}\le 1-\frac{C(x^4+1)}{\sqrt{n}}, \forall x\in B(u),
		\end{equation}
		when $n$ is large enough.
		
		Therefore,
		\begin{equation}\label{bound for h'1}
		|h_1'(x)|\le \begin{cases}
		(\ln n)^2, & x\in B(u)^c;\\
		(1+C)(|x|+1), &x\in B(u).
		\end{cases}
		\end{equation}
		Moreover, 		since $1>\phi(x)-r(x)\ge \frac{1}{n}, \forall x\in B(u)$, by the definition of $h_1$, we have
		\begin{align}
		|h_1(x)|&\le \begin{cases}
		0, &|x|>u\sqrt{\ln n}+\frac{2}{\ln n};\\
		2\ln n, &|x|<u\sqrt{\ln n}+\frac{2}{\ln n}.\label{h1 good bound}
		\end{cases}
		\end{align}
		Due to the fact that for large $n$,
		\begin{equation}\label{ln 2}
		\Big|\ln\Big(1-\frac{C(x^4+1)}{n^{1/2}}-\frac{\sqrt{2\pi}}{n^{(2-u^2)/2}}\Big)\Big|\le \ln 2,\;\;\forall x\in B(u),
		\end{equation}
		then by \eqref{bound for h'1}, we have
		\begin{equation}\label{h_1<C_3}
				\begin{aligned}
		&|\E h_1(G)|\\
		&\le\int_{B(u)} \phi(x) |\ln \phi(x)|dx +\int_{B(u)} \phi(x)\bigg|\ln\Big(1-\frac{C(x^4+1)}{n^{1/2}}-\frac{\sqrt{2\pi}}{n}e^{x^2/2}\Big)\bigg|dx\\
		&\;\;\;\; +\int_{B(u)^c} \phi(x)|h_1(x)|dx\\
		&\le \int_\R  \phi(x) |\ln \phi(x)|dx+\int_{B(u)} \phi(x)\bigg|\ln\Big(1-\frac{C(x^4+1)}{n^{1/2}}-\frac{\sqrt{2\pi}}{n}e^{x^2/2}\Big)\bigg|dx+8\\
		&\le \frac{1}{2}\ln(2\pi e)+\int_{B(u)} \phi(x)\bigg|\ln\Big(1-\frac{C(x^4+1)}{n^{1/2}}-\frac{\sqrt{2\pi}}{n^{1-u^2/2}}\Big)\bigg|dx+8\\
		&\le \frac{1}{2}\ln(2\pi e^2)+\ln 2+8\\
		&:=C_3<\infty,
		\end{aligned}
		\end{equation}
		We next need the following proposition to bound the third term in \eqref{bound for fh''}.
		\begin{proposition}\label{third term}
			\begin{equation}
			\bigg|\frac{1}{\Phi(-|x|)}\int_{-\infty}^{-|x|} h_1(t)\phi(t)dt\bigg|\le C_4 x^21_{|x|\le 2u\sqrt{\ln n}},\;\;\forall x\in\R.
			\end{equation}
		\end{proposition}
		\pf
		We only need to prove the result when $x\in (-u\sqrt{\ln n}-\frac{2}{\ln n}, 0].$  Denote $-u\sqrt{\ln n}-\frac{2}{\ln n}=x_0,$ and $-u\sqrt{\ln n}=x_1.$
		\begin{enumerate}
			\item When $x\in (-u\sqrt{\ln n}-\frac{2}{\ln n}, -u\sqrt{\ln n})$, $h_1$ is linear and it is easy to check that
			\begin{equation}
			\bigg|\int_{-\infty}^{-|x|} h_1(t)\phi(t)dt\bigg|=\bigg|\int_{x_0}^{x} h_1(t)\phi(t)dt\bigg|\le 4\phi(x).
			\end{equation}
			
			\item When $x\in [-u\sqrt{\ln n}, 0]$,
			\begin{equation}
			\bigg|\int_{-\infty}^{-|x|} h_1(t)\phi(t)dt\bigg|=\bigg|\int_{x_0}^{x} h_1(t)\phi(t)dt\bigg|\le 4\phi(x)+\bigg|\int_{x_1}^{x} h_1(t)\phi(t)dt\bigg|.
			\end{equation}
			
			Using \eqref{bound for h'1}, we have
			\begin{align*}
			\bigg|\int_{x_1}^{x} h_1(t)\phi(t)dt\bigg|&\le \int_{x_1}^x (t^2+C_4')\phi(t)dt\\
			&=-t\phi(t)\bigg|_{x_1}^x +\int_{x_1}^{x}  \phi(t)dt+C_4'\Phi(x)\\
			&\le x_1\phi(x_1)+|x|\phi(x)+\P(G\le x)+C_4'\Phi(x)\\
			&\le |x|\phi(x)+\Phi(x)+C_4'\Phi(x).
			\end{align*}
		\end{enumerate}		
		Since $\Phi(-|x|)\ge \frac{1}{|x|}\phi(x)$, we complete the proof.
		\e

		Thus for large $n$, by using  \eqref{bound for fh''}, \eqref{h_1<C_3} and Proposition \ref{third term}, an upper bound is obtained by
		\begin{align}
		&\E f_{h_1}''^2(W_n)\\
		&\le 8\bigg( (1+C)^2\E (|W_n|+1)^2+(\ln n)^4\P(W_n\in B^c(u))\\
		&\;\;\;\;+C_0^2C_3^2\E W_n^2+C_0^2C^2_4\E [W_n^61_{|W_n|\le 2u\sqrt{\ln n}}]+  \E |W_n|^2|h_1(W_n)|^2 \bigg).
		\end{align}

		%		%&\le 2\left( (1+C)^2\E (|W_n|+1)^2+(\ln n)^4\P(W_n\in B(u))+C_0^2\E W_n^2\|h_1\|_\infty \right)\\
		%		&\le C_5(\ln n)^2+8(\ln n)^4(\int_{B^c(u)}p_n(x)-\phi(x)+\phi(x) dx)\\
		%		&\;\;\;\;+8\int_{u\sqrt{\ln n}<|x|<k\ln n} (k\ln n)^2x^2p_n(x)dx\\
		%		&\;\;\;\;
		%		+8\int_{|x|\le u\sqrt{\ln n}} x^2(\ln(\phi(x)-r(x)))^2p_n(x)dx\\
		%		&\le  C_5(\ln n)^2+8(\ln n)^4(\|p_n-\phi\|_1+\P(G\in B^c(u)))\\
		%		&\;\;\;\;+8\frac{C_6}{\sqrt{n}}\int_{u\sqrt{\ln n}<|x|<k\ln n} (k\ln n)^2x^2dx\\
		%		&\;\;\;\;+8\int_{|x|\le u\sqrt{\ln n}} x^2(\ln(\phi(x)-r(x)))^2|p_n(x)-\phi(x)|dx\\
		%		&\;\;\;\;+8\int_{|x|\le u\sqrt{\ln n}}x^2(\ln(\phi(x)-r(x)))^2\phi(x)dx\\
		%		&\le C_5(\ln n)^2+C_7\frac{(\ln n)^5}{\sqrt{n}}+8\int_{B(u)} x^2\phi(x)|\ln \phi(x)|^2 dx\\
		%		 &\;\;\;\;+8\int_{B(u)}x^2\phi|\ln(1-\frac{C(x^m+1)}{n^{1/2}}-\frac{\sqrt{2\pi}}{n^{(k-u^2)/2}})|^2dx\\
		%		&\le C_5(\ln n)^2+C_7\frac{(\ln n)^5}{\sqrt{n}}+16\int_\R x^2(x^4+\ln(\sqrt{2\pi})^2) \phi(x) dx +8(\ln 2)^2\E(G^2)\\
		%		&\le C_8(\ln n)^2. \label{f''}
		%%&        \le C_3(\ln n)^2+2(\ln n)^4(\int_{B(u)}p_n-\phi+\phi)\\
		%%&\le C_3(\ln n)^2+2(\ln n)^4(\|p_n-\phi\|_1+\P(G\in B(u)))\\
		%%&\le C_3(\ln n)^2+C_4\frac{k^5(\ln n)^5}{\sqrt{n}}\\
		%%&\le C_5(\ln n)^2.
		%		\end{align}
		Now we bound the terms at the RHS of the above inequality as follows.
		\begin{enumerate}
			\item Clearly, $(1+C)^2\E (|W_n|+1)^2+C_0^2C_3^2\E W_n^2$ is bounded by a constant.
			
			\item Next, since
			\begin{equation}
			C_0^2C^2_4\E [W_n^61_{|W_n|\le 2u\sqrt{\ln n}}]\le C_0^2C^2_44u^4(\ln n)^2\E [W_n^2]=C_0^2C^2_44u^4(\ln n)^2,
			\end{equation}
			we obtain that
			\begin{align}
			&8\left( (1+C)^2\E (|W_n|+1)^2+C_0^2C_3^2\E W_n^2+ _0^2C^2_4\E [W_n^61_{|W_n|\le 2u\sqrt{\ln n}}]\right)\\
			& \le C_5(\ln n)^2.
			\end{align}
			
			\item Then, we have the following bound
			\begin{align}
			\P(W_n\in B^c(u))&=\int_{B^c(u)}p_n(x)-\phi(x)+\phi(x) dx\\
			&\le \|p_n-\phi\|_1+\P(G\in B^c(u)).
			\end{align}
			Using	the fact that $\|p_n-\phi\|_1=O(\ln n/\sqrt{n})$, we obtain
			\begin{equation}
			8(\ln n)^4\P(W_n\in B^c(u))\le C_6\frac{(\ln n)^5}{\sqrt{n}}.
			\end{equation}
			
			\item Finally, since $u\sqrt{\ln n}+\frac{2}{\ln n}\le u\sqrt{\ln n}+2\ln n$ when $n$ is large, we have
			\begin{align}
			&\E |W_n|^2|h_1(W_n)|^2=\int_{|x|<2\ln n} h^2_1(x)x^2p_n(x)dx\\
			&\le \int_{u\sqrt{\ln n}<|x|<2\ln n} (2\ln n)^2x^2p_n(x)dx\\
			&\;\;\;\;
			+\int_{|x|\le u\sqrt{\ln n}} x^2(\ln(\phi(x)-r(x)))^2p_n(x)dx.
			\end{align}
			On the one hand, by \eqref{local bounds} we have $p_n(x)\le \frac{C_7}{\sqrt{n}}$, $\forall x\text{ } s.t.\text{ }u\sqrt{\ln n}<|x|<2\ln n$ when $n$ is large. Thus
			\begin{equation}
			\int_{u\sqrt{\ln n}<|x|<2\ln n} (2\ln n)^2x^2p_n(x)dx\le \frac{C_7(\ln n)^5}{\sqrt{n}}.
			\end{equation}
			
			On the other hand,
			\begin{align}
			&\int_{|x|\le u\sqrt{\ln n}} x^2(\ln(\phi(x)-r(x)))^2p_n(x)dx\\
			&\le \int_{|x|\le u\sqrt{\ln n}} x^2(\ln(\phi(x)-r(x)))^2|p_n(x)-\phi(x)|dx\\
			&\;\;\;\;+
			\int_{|x|\le u\sqrt{\ln n}}x^2(\ln(\phi(x)-r(x)))^2\phi(x)dx\\
			&\le \int_{B(u)}x^2\phi(x)|\ln\left( 1-\frac{C(x^m+1)}{n^{1/2}}-\frac{\sqrt{2\pi}}{n^{(2-u^2)/2}}\right) |^2dx\\
			&\;\;\;\;+\int_{B(u)} x^2\phi(x)|\ln \phi(x)|^2 dx\\
			&\le  (\ln 2)^2\E(G^2)+2\int_\R x^2(x^4+\ln(\sqrt{2\pi})^2) \phi(x) dx,
			\end{align}
			where we use \eqref{ln 2}   in the last inequality. Therefore,
			\begin{equation}
			8\E |W_n|^2|h_1(W_n)|^2\le C_8.
			\end{equation}
		\end{enumerate}
		In conclusion, we obtain the following upper bound
		\begin{equation}\label{f''}
		\E f_{h_1}''^2(W_n)\le C_9(\ln n)^2.
		\end{equation}

		As for $f'''_{h_1}$, first we have
		$\|h'_1\|_\infty\le (\ln n)^2.$
		Recall that in \eqref{1+C},    when $n$ is large we have,
		$$0<\frac{(1-C(x^m+1)   /\sqrt{n})\phi(x)}{\phi(x)-r(x)}\le 1+C,\;\;\forall x\in B(u).$$
		By some appropriate smoothing, we can obtain
		\begin{align*}
		|h''_1(x)|\le C_{10}(x^2+1),\;\;\forall x\in B(u),
		\end{align*}
		and
		%$$|h''_1(x)|\le n^{1/4},\;\; h_1'(x)h_1''(x)\le 0,\;\;\forall x\in B(u)^c$$
		$$|h''_1(x)|\le n^{1/4},\;\;\forall x\in B(u)^c$$
		without change any the above bounds. In fact, by this bound for $h_1''$, we only need to smooth $h_1$ at small intervals
		$$(-u\sqrt{\ln n}-\frac{(1+C)(u\sqrt{\ln n}+1)+(\ln n)^2}{n^{1/4}}, -u\sqrt{\ln n}),$$
		$$(u\sqrt{\ln n}, u\sqrt{\ln n}+\frac{(1+C)(u\sqrt{\ln n}+1)+(\ln n)^2}{n^{1/4}}),$$
		and the neighborhoods of two zero arriving points in the outer area $|x|>u\sqrt{\ln n}$, and on these intervals, $h_1(x)$ can only increase $\varepsilon$ that satisfies $\varepsilon\le \|h_1'\|_\infty\frac{(1+C)(u\sqrt{\ln n}+1)+(\ln n)^2}{n^{1/4}}\le 1$ when $n$ is large, \textit{i.e.} the bound \eqref{h1 good bound} still holds.

		Therefore, by \eqref{Stein f'''} and \eqref{bound for fh'''}
		\begin{equation}\label{I_2,2}
		I_2\le C_{11}\frac{\ln n}{\sqrt{n}}+C_{12}\frac{(\ln n)^{5/2}+n^{1/4}}{n}.
		\end{equation}

		\item Upper bound of $I_3$. According to \eqref{h1 good bound} and the fact that for large $n$, $p_n(x)+\phi(x)\le \frac{C_{13}}{\sqrt{n}}, \forall x\in B(u)^c$,
		\begin{align}
		I_3&\le \int_{u\sqrt{\ln n}<|x|<u\sqrt{\ln n}+2/\ln n} 2\ln n|p_n(x)-\phi(x)|dx\\
		&\le \int_{u\sqrt{\ln n}<|x|<u\sqrt{\ln n}+2/\ln n} 2\ln n(p_n(x)+\phi(x))dx\\
		&\le \frac{C_{13}}{\sqrt{n}}\int_{u\sqrt{\ln n}<|x|<u\sqrt{\ln n}+2/\ln n} 2\ln ndx\\
		&\le \frac{2^2C_{13}}{\sqrt{n}}.\label{I_3,2}
		\end{align}
		
		\item Upper bound of $I_4$.  By Proposition \ref{Thm has Pro} and $u\ge 1$, we have
		\begin{align}
		I_4&\le \int_{B(u)^c\cap A^c}(p_n(x)-\phi(x)) \ln p_n(x) dx\\
		&\le \int_{B(u)^c\cap A^c}(\phi(x)-p_n(x))\left( |\ln k_1|+ n^sk_2e^{\frac{x^2}{2a} }\right)dx\\
		&\le \int_{B(u)^c}\phi(x)\left( |\ln k_1|+ n^sk_2e^{\frac{x^2}{2a} }\right)dx\\
		&\le k_3 \int_{B(u)^c} n^se^{-\frac{x^2}{2}(1-\frac{1}{a}) }dx\\
		&\le k_4 \int_{B(u\sqrt{1-\frac{1}{a}})^c} n^se^{-\frac{x^2}{2} }dx\\
		&\le k_5 n^{s-u^2(1-1/a)},\label{I_4,2}
		\end{align}
	\end{enumerate}
	where $k_3,k_4,k_5$ are constants.
	Recall that $2>\frac{a}{a-1}(\frac{1}{2}+s)$ according to Property \ref{e^e assumption}, thus we can find $u\in [1,\sqrt{2})$ \textit{s.t.}
		\begin{equation}
		s-u^2(1-1/a)\le -\frac{1}{2}.
		\end{equation}

	Combining \eqref{I_1,2}, \eqref{I_2,2}, \eqref{I_3,2} and  \eqref{I_4,2}, we obtain the Theorem \ref{bound for d under Stein}.
	
	\e
	
%	\begin{remark}
%		If the conditions in Proposition \ref{full density assumption} hold, \textit{i.e.} the Assumption \ref{part A assume} holds with $A_n=\{1,2,\cdots, n\}$, then any $k\ge 2$  meets the conditions in Theorem \ref{bound for d under Stein} and we can choose $u=\sqrt{\frac{11}{6}}$ in our proof.
%	\end{remark}

	\section{Some Remarks on Rate Improvement}\label{section4}
	In this section, we discuss some methods to improve the convergence rate $O(\ln n/\sqrt{n})$ in Theorem \ref{bound for d under Stein}.

	\subsection{Assumptions for $O(\frac{1}{\sqrt{n}})$ Convergence Rate}
	\textbf{If we have $\sup_n \E W_n^6<\infty$}, then \eqref{f''}  can be bounded  by a constant not depending on $n$ and thus \eqref{I_2,2} can be modified as $O(\frac{1}{\sqrt{n}})$. Then
	\begin{equation}
	d(W_n)=O\left(\frac{1}{\sqrt{n}}\right).
	\end{equation}
	
	Furthermore, \textbf{if $\V(W_n-\frac{X_i}{\sum_{j=1}^n \V(X_j)})$ tends to $1$ uniformly on $i$}, we can  also obtain the $O(\frac{1}{\sqrt{n}})$ convergence rate. Indeed, we modify the upper bound of \eqref{3.2}  with $h=h_1$ and use \eqref{bound for fh''}, Proposition  \ref{third term} and \cite[Lemma 2.8]{steinbook} to obtain the following estimation.
	\begin{equation}
	\begin{aligned}
	&|\E h_1(W_n)- \E h_1(G)|\le \E[f_{h_1}''(W_n)\Delta]+\E[f_{h_1}'''(\beta)\Delta^2]\\
	&=\E[f_{h_1}''(W_n-\xi_I)\Delta]+\E[f_{h_1}'''(\hat{\beta})\Delta]+\E[f_{h_1}'''(\beta)\Delta^2]\\
	&\le  \sum_{i=1}^n (\E |\xi_i|^3)(\E|h'_1(W_n-\xi_i)|+\E| h_1(G)|\E|W_n-\xi_i|+\E|(W_n-\xi_i)h_1'(W_n-\xi_i)|+\sqrt{\E W_n^4})\\
	&\;\;\;\;+ \|f_{h_1}'''\|_\infty \sum_{i=1}^n \E \xi_i^2 \E|(\xi_i-\xi_i^*)\xi_i|+ \gamma\|f_{h_1}'''\|_\infty \frac{1}{n},
	\end{aligned}
	\end{equation}
	where $\hat{\beta}$ is a random variable between $W_n-\xi_I$ and $W_n$; $\beta$ is a random variable between $W_n$ and $W_n^*$; $\xi_i, \xi_i^*, I$ and $\gamma$ are the same in Lemma \ref{bound of delta^2}.
	$h_1$ does not depend on the density function $p_n$, thus we can treat $\frac{W_n-\xi_i}{\V(W_n-\xi_i)}$ as a sum of $n-1$ independent random variables and then  we can follow the previous discussion  for the bound of $I_2$ to obtain
	\begin{equation}
	\begin{aligned}
	&\sum_{i=1}^n (\E |\xi_i|^3)(\E|h'_1(W_n-\xi_i)|+\E| h_1(G)|\E|W_n-\xi_i|+\E|(W_n-\xi_i)h_1'(W_n-\xi_i)|+\sqrt{\E W_n^4})\\
	&\le \sum_{i=1}^n (\E |\xi_i|^3)(C_{13}+\sqrt{\E W_n^4})\\
	&\le \sum_{i=1}^n  ( \E |\xi_i|^3)\left( C_{13}+\sqrt{\frac{\sum_{j=1}^n \E X_j^4+3(\sum_{j=1}^n \E X_j^2)^2}{(\sum_{j=1}^n \E X_j^2)^2}}\right) \\
	&\le \sum_{i=1}^n (\E |\xi_i|^3)\left( C_{13}+C_{14}\sqrt{\frac{ n}{(\sum_{j=1}^n \E X_j^2)^2}+1}\right) \\
	&=\frac{\sum_{i=1}^n \E |X_i|^3}{(\sum_{i=1}^n \E |X_i|^2)^{3/2}}\left( C_{13}+C_{14}\sqrt{\frac{ n}{(\sum_{i=1}^n \E X_i^2)^2}+1}\right) ,
	\end{aligned}
	\end{equation}
	where we use the assumption that $\V(W_n-\xi_i)$ is uniformly close to 1 about $i$, when $n$ is large.  By the definition of the zero bias, $\E (\xi_i^*)^2=\frac{1}{3\E\xi_i^2}  \E \xi_i^4 $, then we conclude by the condition $\sup_i \E X_i^4<\infty$.
	\begin{equation}
	\begin{aligned}
	&|\E h_1(W_n)- \E h_1(G)|\\
	&\le \frac{n}{(\sum_{i=1}^n \E |X_i|^2)^{3/2}}\left( C_{13}+C_{14}\sqrt{\frac{ n}{(\sum_{i=1}^n \E X_i^2)^2}+1}\right) +C_{15}\frac{n^{5/4} }{(\sum_{i=1}^n \E X_i^2)^2}+C_{16}n^{-3/4}.
	\end{aligned}
	\end{equation}
	Particularly, when $X_i$s  have equal variances, then $|\E h_1(W_n)- \E h_1(G)|=O(1/\sqrt{n})$ and thus $d(W_n)=O(1/\sqrt{n})$.

	\subsection{Without Property  \ref{e^e assumption}}
	
	If we delete   Property \ref{e^e assumption} (or Assumption \ref{part A assume}), we can still have
	\begin{equation}
	d(W_n)=O\left(\frac{\ln n}{\sqrt{n}}\right)+I_4,
	\end{equation}
	since $I_1, I_2,$ and $I_3$ are not related to Property \ref{e^e assumption} (or Assumption \ref{part A assume}), where $I_1, I_2, I_3$ and $I_4$ are defined in Section \ref{section 3.2}. For the convergence of $I_4$, we can only make an assumption:
	$T_n(x):=\phi(x)\ln\frac{\phi(x)}{p_n(x)}$ are uniformly integrable about $n$. Recall that $\int_\R T_n(x)dx=\KL(G\|W_n)$, the above assumption is natural since we always need $T_n$ to be integrable. However, the above assumption can not deduce  convergence rate of $I_4$, thus also for $d(W_n)$.
	
	\subsection{Discussion for $O(\frac{1}{n})$ Convergence Rate} 
	If we consider function $\hat{g}(x):=\ln(p_n(x)/\phi(x))$ instead of $g(x)=\ln p_n(x)$,  then by the analogous method with same truncation, we can obtain
	\begin{equation}\label{dn as 1/n with tial}
	d(W_n)= O\left(\frac{1}{n}\right)+\int_{|x|>u\sqrt{\ln n}}p_n(x)\ln \phi(x)dx.
	\end{equation}
	In fact, under the conditions in Theorem \ref{bound for d under Stein}, the convergence rate of  KL-divergence
	\begin{align}
	&\KL(W_n\|G)\\
	&=\int_{|x|>u\sqrt{\ln n}}p_n(x)\ln p_n(x)dx-\int_{|x|>u\sqrt{\ln n}}p_n(x)\ln \phi(x)dx\\
	&\;\;\;\;+\int_{|x|<u\sqrt{\ln n}}p_n(x)\ln(p_n(x)/\phi(x))dx\\
	&:=\hat{R}_1-\hat{R}_2+\hat{R}_3
	\end{align}
	is $\frac{1}{n}$ since the results in \cite{Esseenbound} and the following  well known inequality in \cite{S-EPI}.
	\begin{lemma}\label{control D by VJ}
		Given $\xi$, $\E \xi=0, \V(\xi)<\infty$, then
		\begin{equation}
		D(\xi)\le \frac{1}{2}\ln(\V(\xi) J(\xi)).
		\end{equation}
	\end{lemma}
	Additionally,  $\hat{R}_1$ can be estimate by Property \ref{e^e assumption} and $\hat{R}_3$ may be almost canceled by the term
	$\int_{|x|\le u\sqrt{\ln n}}\phi(x)\ln(p_n(x)/\phi(x))dx$ in $\KL(G\|W_n)$.  Note that the remaining term $\hat{R_2}$ has the same order of
	\begin{equation}
	\hat{R}(n):=\int_{|x|>u\sqrt{\ln n}}x^2p_n(x)dx,
	\end{equation}
	thus we have
	\begin{equation}
	d(W_n)= O\left(\frac{1}{n}\right)+\hat{R}(n).
	\end{equation}
	Unfortunately, the estimation for tail variance $\hat{R}(n)$ can not be more precise with known tools, such as Chebyshev's inequality. But combining our main result Theorem \ref{bound for d under Stein} and \eqref{dn as 1/n with tial}, we can deduce an new estimation for tail variance by: if the conditions in Theorem \ref{bound for d under Stein} hold, then
	\begin{equation}
	\int_{|x|>u\sqrt{\ln n}}x^2p_n(x)dx=O\left( \frac{\ln n}{\sqrt{n}}\right) ,\;\;\forall u>0.
	\end{equation}

	On the other hand, we can add stronger moment condition to overcome the unknown speed of $\hat{R}(n)$.  Particularly, if $\exists \beta>0, \sup_i \E(\exp(\frac{X_i^2}{\beta}))<\infty$, then it is easy to check that
	\begin{equation}\label{order of d as 1/n}
	d(W_n)=O\left( \frac{1}{n}\right) ,\;\; n\to\infty.
	\end{equation}
	Since local limit theorem says that under these conditions $\|p_n-\phi\|_1\sim \frac{1}{\sqrt{n}}$, \textit{i.e.} has the same order with $ \frac{1}{\sqrt{n}}$, if $R_1(x)\neq 0$, therefore, thanks to Pinsker's inequality \eqref{Pinsker}, we conclude that \eqref{order of d as 1/n} implies 'the best order estimation' of the rate theorem of symmetric KL-divergence.
	%				\begin{equation}
	%				d(W_n)\sim \frac{1}{n};\;\;  \KL(W_n\|G)\sim \frac{1}{n};\;\; \KL(G\|W_n)\sim \frac{1}{n}.
	%				\end{equation}
	%				
	
	\section*{Acknowledgments}
	\noindent We thank Zhi-Ming Ma  for his careful reading of proofs in the manuscript. We thank Qi-Man Shao for his useful comments on this article. We thank the reviewers for their  helpful suggestions.

	\section*{Funding}
	\noindent Liu S.-H. was partially supported by National Nature Science Foundation of China NSFC 12301182.
	
	\noindent Yao L.-Q. was partially supported by National Key R\&D Program of China No. 2023YFA1009603.

	\appendix
	\setcounter{proposition}{0}
	\renewcommand{\theproposition}{A.\arabic{proposition}}
	
	\section{Assumption \ref{part A assume} implies Property \ref{e^e assumption}}\label{App1}
	Here we   prove that under the Assumption \ref{part A assume}, Property \ref{e^e assumption} holds for $W_n$ immediately. 
	Firstly, we consider a more special case of  Assumption \ref{part A assume}.
	\begin{proposition}\label{full density assumption}
		Given a   random variable sequence $X_1, X_2, \cdots, X_n$. %\textit{s.t.} $\E X_i=0, \sup_i J(X_i)<J<\infty, \sup_i \E|X_i|^{4+\delta_0}<M<\infty$ for some $k\in\mathbb{N}^+, \delta_0>0, \forall i=1,2,\cdots,n$.
		If   the  density functions $\rho_1, \rho_2, \cdots, \rho_n$ of $X_1, X_2, \cdots, X_n$ satisfy
		\begin{equation}
		\rho_i(x)\ge l_1e^{-l_2\frac{x^2}{2}},\forall x\in \R, \forall i=1,2,\cdots,n,
		\end{equation}
		for some $l_1>0, l_2>0$, then $\forall v_0>0, a_0>1$, \eqref{p>e} holds with
			\begin{equation}
			a=a_0, \;k_1=l_1, \;k_2=l_2+l_3, \;s=1, \;v=v_0
			\end{equation} 	
		for the density function $q_n$ of
		\begin{equation}
		V_n:=\frac{1}{\sqrt{n}}\sum_{i=1}^n X_i,
		\end{equation}
		when   $n$ is large enough (depends only on $a_0$).
	\end{proposition}
	\pf Note that the density of $(X_1+X_2)/\sqrt{2}$ satisfies that
	\begin{align*}
	p_{\frac{X_1+X_2}{\sqrt{2}}}(x)&= \int_\R \sqrt{2}\rho_1(\sqrt{2}x-t)\rho_2(t)dt\\
	&\ge \sqrt{2}\int_\R l_1e^{-l_2\frac{(\sqrt{2}x-t)^2}{2}}l_1e^{-l_2\frac{t^2}{2}}dt\\
	&= l^2_1\sqrt{\frac{2\pi}{l_2}}e^{-l_2\frac{x^2}{2}}
	\end{align*}
	and then we prove by induction: when $q_{n}(x)\ge l_1\left( l_1\sqrt{\frac{2\pi}{l_2}}\right)^{n-1}e^{-l_2\frac{x^2}{2}}$, we can deduce that
	\begin{equation*}
\begin{aligned}
q_{n+1}(x)&=\int_\R \frac{\sqrt{n+1}}{\sqrt{n}}q_n\left( \frac{\sqrt{n+1}}{\sqrt{n}}(x-t)\right) \sqrt{n+1}\rho_{n+1}(\sqrt{n+1}t)dt\\
	&\ge \int_\R \frac{\sqrt{n+1}}{\sqrt{n}}l_1\left( l_1\sqrt{\frac{2\pi}{l_2}}\right)^{n-1}e^{-l_2\frac{\left( \frac{\sqrt{n+1}}{\sqrt{n}}(x-t)\right)^2}{2}} \sqrt{n+1}l_1e^{-l_2\frac{(n+1)t^2}{2}}dt\\
		&=\left( l_1\sqrt{\frac{2\pi}{l_2}}\right)^{n+1}\int_\R \frac{\sqrt{n+1}}{\sqrt{n}} \sqrt{\frac{l_2}{2\pi}}e^{-l_2\frac{\left( \frac{\sqrt{n+1}}{\sqrt{n}}(x-t)\right)^2}{2}}\sqrt{n+1}\sqrt{\frac{l_2}{2\pi}}e^{-l_2\frac{(n+1)t^2}{2}}dt\\
		&=\left( l_1\sqrt{\frac{2\pi}{l_2}}\right)^{n+1} \sqrt{\frac{l_2}{2\pi}}e^{-l_2\frac{x^2}{2}}\\
		&=l_1\left( l_1\sqrt{\frac{2\pi}{l_2}}\right)^{n} e^{-l_2\frac{x^2}{2}}.
\end{aligned}
	\end{equation*}
	Therefore,
	 we conclude that  the density $q_n$ of $V_n$ satisfies 
	\begin{align}
	q_n(x)&\ge l_1\left( l_1\sqrt{\frac{2\pi}{l_2}}\right)^{n-1}e^{-l_2\frac{x^2}{2}}\ge l_1e^{-nl_3}e^{-l_2\frac{x^2}{2}}\\
	&\ge l_1e^{-(l_2+l_3)ne^{x^2/2a_0},\;\;\forall |x|>v_0\sqrt{\ln n}}, %|x|>\sqrt{\ln n},
	\end{align}
	once $n$ is large enough, where $l_3:=\max\{1, -\ln\left( l_1\sqrt{\frac{2\pi}{l_2}} \right)  \}.$
	Therefore, $q_n$ satisfies \eqref{p>e} for $k_1=l_1, s=1, k_2=l_2+l_3, a=a_0, v=v_0$.
	\e

	Proposition \ref{full density assumption} shows an extreme case of Assumption \ref{part A assume} when $A=\{1,2,\cdots, n\}$. Based on it, we conclude the following result.
	\begin{proposition}\label{Thm has Pro}
		%	if $W_n=\frac{1}{\sqrt{n}}\sum_{i=1}^n X_i$, $\V(X_i)=1, \forall i,$ and $X_i$s satisfy the Assumption \ref{part A assume}, then $W_n$ satisfies Property \ref{e^e assumption}.
		%	Given $W_n=\frac{\sum_{i=1}^n X_i}{\sqrt{\sum_{j=1}^n \V(X_j)}}$,
		%	where $\E X_i=0, \sup_i J(X_i)<J<\infty, \sup_i \E|X_i|^{4+\delta_0}<M<\infty$ for some $k\in\mathbb{N}^+, \delta_0>0, \forall i=1,2,\cdots,n$.
		%	
		%	If $n$ is large enough and $\{X_1, X_2, \cdots, X_n\}$ satisfies the Assumption \ref{part A assume}, then
		Under the conditions in Theorem \ref{bound for d under Stein},
		$W_n$  satisfies the  Property \ref{e^e assumption} with
		\begin{equation}
			k_1=\sqrt{\frac{1}{4J}} l_1,\; s=1,\; k_2=l_2+l_3,\; a=\frac{\delta_1}{M^{\frac{2}{4+\delta_0}}},\; v=1.
			\end{equation}
		when $n$ is sufficient large.
	\end{proposition}
	\pf   Let $Y_1=\frac{1}{\sqrt{|A_n|}}\sum_{i\in A_n} X_i$, $Y_2=\frac{1}{\sqrt{n-|A_n|}}\sum_{i\notin A_n} X_i$, where $A_n$ is defined in Assumption \ref{part A assume}, then
	$$W_n=\frac{\sqrt{|A_n|}}{B_n}Y_1+\frac{\sqrt{n-|A_n|}}{B_n}Y_2,$$
	where $B_n=\sqrt{\sum_{j=1}^n \V(X_j)}.$
	Denote the densities of $Y_1, Y_2$ by $i_1(x)$ and $i_2(x)$ respectively. By  Proposition \ref{full density assumption} with $a_0=4$, it is easy to check that $\forall v'>0$,  if  $n$ is large enough, 
	\begin{equation}\label{2.1_A_n}
	i_1(x)\ge l_1e^{-(l_2+l_3)|A_n|e^{x^2/8}},\;\;\forall |x|>v'\sqrt{\ln |A_n|},
	\end{equation}
	where $l_3:=\max\left\{1, -\ln\left( l_1\sqrt{\frac{2\pi}{l_2}}\right) \right\}.$
	
	Given $v=1, r=1.5$ and $n$ is large enough, for $x$ \textit{s.t} $|x|>v\sqrt{\ln n}$,
		\begin{align}
		p_n(x)&=\int_{t\in\R} \frac{B_n}{\sqrt{|A_n|}} i_1\left( \frac{B_n}{\sqrt{|A_n|}}(x-t)\right) \cdot \frac{B_n}{\sqrt{n-|A_n|}} i_2\left( \frac{B_n}{\sqrt{n-|A_n|}}t\right) dt\\
		&\ge \int_{|t|\le (r-1)v\sqrt{\ln n}}\frac{B_n}{\sqrt{|A_n|}} i_1\left( \frac{B_n}{\sqrt{|A_n|}}(x-t)\right) \cdot \frac{B_n}{\sqrt{n-|A_n|}} i_2\left( \frac{B_n}{\sqrt{n-|A_n|}}t\right) dt\\
		&\overset{(a)}{\ge} \int_{|t|\le (r-1)v\sqrt{\ln n}} \frac{B_n}{\sqrt{|A_n|}} l_1e^{-(l_2+l_3)|A_n|e^{\frac{B^2_n}{8|A_n|}(x-t)^2}} \cdot \frac{B_n}{\sqrt{n-|A_n|}} i_2\left( \frac{B_n}{\sqrt{n-|A_n|}}t\right) dt\\
		&\overset{(b)}{\ge} \int_{|t|\le (r-1)v\sqrt{\ln n}} \frac{B_n}{\sqrt{|A_n|}} l_1e^{-(l_2+l_3)ne^{\frac{B^2_n}{8\delta_1 n}(rx)^2}} \cdot \frac{B_n}{\sqrt{n-|A_n|}} i_2\left( \frac{B_n}{\sqrt{n-|A_n|}}t\right) dt\\
		&\ge \left( \frac{B_n}{\sqrt{|A_n|}} l_1e^{-(l_2+l_3) ne^{\frac{r^2B^2_nx^2}{8\delta_1 n}}}  \right) \int_{|t|\le (r-1)v\sqrt{\ln n}/2}  \frac{B_n}{\sqrt{n-|A_n|}} i_2\left( \frac{B_n}{\sqrt{n-|A_n|}}t\right) dt\\
		&=\left( \frac{B_n}{\sqrt{|A_n|}} l_1e^{-(l_2+l_3)ne^{\frac{r^2B^2_nx^2}{8\delta_1 n}}}  \right) \P\left( |Y_2|\le \frac{B_n}{\sqrt{n-|A_n|}}\frac{(r-1)v\sqrt{\ln n}}{2}\right)\\
		&\overset{(c)}{\ge}  \sqrt{\frac{1}{J}} l_1e^{-(l_2+l_3)ne^{\frac{ r^2M^{\frac{2}{4+\delta_0}}}{4\delta_1 }(x^2/2)}}   \P\left( |Y_2|\le \sqrt{\frac{1}{J(1-\delta_1)}}\frac{(r-1)v\sqrt{\ln n}}{2}\right)\\
		&\overset{(d)}{\ge}   \sqrt{\frac{1}{J}} l_1e^{-(l_2+l_3)ne^{\frac{ M^{\frac{2}{4+\delta_0}}}{\delta_1 }(x^2/2)}}   \P\left( |Y_2|\le \sqrt{\frac{1}{J(1-\delta_1)}}\frac{(r-1)v\sqrt{\ln n}}{2}\right),
		\end{align}
		where $(a)$ holds according to \eqref{2.1_A_n}, CR bound $\V(\xi)J(\xi)\ge 1$ for any random variable,  and the fact that when $|x|>v\sqrt{\ln n}, |t|\le (r-1)v\sqrt{\ln n}$,
		\begin{align}
		\frac{B_n}{\sqrt{|A_n|}}|x-t|&\ge \frac{\sqrt{\sum_{i=1}^n \V(X_i)}}{\sqrt{n}}(2-r)v\sqrt{\ln n}\\
		&\ge \frac{\sqrt{\sum_{i=1}^n 1/J(X_i)}}{\sqrt{n}}(2-r)v\sqrt{\ln n}\\
		&\ge \frac{(2-r)v}{\sqrt{J}}\sqrt{\ln n}.
		\end{align}
		$(b)$ holds by the fact $\delta_1 n\le |A_n|\le n$ and $(x-t)^2\le (rx)^2$,
	and $(c)$ holds since   $\sup_i J(X_i)<J<\infty, \sup_i \E|X_i|^{4+\delta_0}<M<\infty$ and the fact that for any random variable $\xi$,
	\begin{equation}
	\frac{1}{J(\xi)}\le \V(\xi)\le (\E|\xi|^m)^{\frac{2}{m}},\;\;\forall m>2.
	\end{equation}
	$(d)$ is deduced by $r<2$.
	When $n$ is large enough,
	\begin{equation}
	\P\left( |Y_2|\le \sqrt{\frac{1}{J(1-\delta_1)}}\frac{(r-1)v\sqrt{\ln n}}{2}\right)>\frac{1}{2}
	\end{equation}
	by Berry-Esseen bound (see \cite[Theorem 3.3]{steinbook} as an example), thus we conclude that
	\begin{align}
	p_n(x)&\ge \sqrt{\frac{1}{J}} l_1e^{-(l_2+l_3)ne^{\frac{ M^{\frac{2}{4+\delta_0}}}{4\delta_1 }(x^2/2)}}   \P\left( |Y_2|\le \sqrt{\frac{1}{J(1-\delta_1)}}\frac{(r-1)v\sqrt{\ln n}}{2}\right)\\
	&\ge \sqrt{\frac{1}{4J}} l_1e^{-(l_2+l_3)ne^{\frac{ M^{\frac{2}{4+\delta_0}}}{4\delta_1 }(x^2/2)}},\;\;\forall |x|>v\sqrt{\ln n}.
	\end{align}
	\textit{i.e.} we can take
	\begin{equation}
	k_1=\sqrt{\frac{1}{4J}} l_1,\; s=1,\; k_2=(l_2+l_3),\; a=\frac{\delta_1}{M^{\frac{2}{4+\delta_0}}},\; v=1.
	\end{equation}
	
	Note that we can alway assume that $M$ is small enough \textit{s.t.}
		\begin{equation}\label{M}
		a=\frac{\delta_1}{M^{\frac{2}{4+\delta_0}}}>1,\;\; 
		\frac{2(a-1)}{a}=2\left( 1-\frac{M^{\frac{2}{4+\delta_0}}}{\delta_1}\right) >\frac{3}{2}=\frac{1}{2}+s,
		\end{equation}
		since if $M$ does not satisfies \eqref{M}, we can find $M_0>0$ satisfying \eqref{M} and change $X_i$ to $\hat{X}_i=\left( \frac{M_0}{M}\right)^{\frac{1}{4+\delta_0}}X_i, i=1,2,\cdots $ with parameters 
		$$(J, M,  \delta_0,\delta_1,l_1,l_2)$$
		changing to 
		$$\left[ \left( \frac{M}{M_0}\right)^{\frac{2}{4+\delta_0}}J, M_0,  \delta_0,\delta_1,\left( \frac{M}{M_0}\right)^{\frac{1}{4+\delta_0}}l_1,\left( \frac{M}{M_0}\right)^{\frac{1}{4+\delta_0}}l_2\right], $$
		and $W_n$ is invariable.
	\e

\end{document}